\documentclass[11pt]{amsart}

\usepackage{amsmath, amssymb}
\usepackage{amsfonts}
\usepackage{graphicx}
\usepackage{mathrsfs}
\usepackage[arrow,matrix,curve,cmtip,ps]{xy}
\usepackage{amsthm}
\usepackage{enumerate}
\usepackage{hyperref}
\usepackage{tikz}

\allowdisplaybreaks

\newtheorem*{theorem*}{Theorem}
\theoremstyle{remark}

\numberwithin{equation}{section}

\newcommand{\g}{\mathcal{G}}
\newcommand{\go}{\mathcal{G}^0}
\newcommand{\gl}{\mathcal{G}^1}
\newcommand{\lrg}{L_R(\mathcal{G})}
\newcommand{\rg}{r_{\mathcal{G}}}
\newcommand{\sg}{s_{\mathcal{G}}}
\newcommand{\qGo}{\Phi(G^0)}
\newcommand{\qgo}{\Phi(\mathcal{G}^0)}
\newcommand{\qgl}{\Phi(\mathcal{G}^1)}
\newcommand{\qsg}{\Phi_{\mathrm{sg}}(G^0)}
\newcommand{\qg}{\mathcal{G}/(H,S)}
\newcommand{\lqg}{L_R(\mathcal{G}/(H,S))}
\newcommand{\I}{I_{(H,S)}}
\newcommand{\el}{\mathcal{EL}_A(R)}

\usepackage[top=3.7cm,right=3.7cm,bottom=3.1cm,left=3.7cm]{geometry}

\begin{document}
\title[The Leavitt Path Algebras of  Ultragraphs]{The Leavitt Path Algebras of Ultragraphs}

\author[M. Imanfar, A. Pourabbas and H. Larki]{M. Imanfar, A. Pourabbas AND H. Larki}
\address{Faculty of Mathematics and Computer Science,
	Amirkabir University of Technology, 424 Hafez Avenue, 15914
	Tehran, Iran.}

\email{m.imanfar@aut.ac.ir, arpabbas@aut.ac.ir}
\address{Department of Mathematics, Faculty of Mathematical Sciences and Computer,
Shahid Chamran University of Ahvaz, Iran}
\email{h.larki@scu.ac.ir}

%\thanks{This research was supported by NSF Grant DMS-1234567}

%\date{\today}



\begin{abstract}
We introduce  the   Leavitt path algebras of ultragraphs  and we  characterize their ideal structures.
	We then use this notion  to introduce and  study the algebraic analogy of Exel-Laca algebras. 
\end{abstract}

\maketitle

%\baselineskip=.45cm

%%%%%%%%%%%%%%%%%%%%%%%%%%%%%%%%%%%%%%%%%%%%%%%%
\section{Introduction}

The Cuntz-Krieger algebras were introduced and studied in \cite{cun} for binary-valued matrices  with finite index. Two immediate and important extensions of the Cuntz-Krieger algebras are: (1) the class of $C^*$-algebras  associated to (directed) graphs  \cite{kum1,kum2,bat2,fow} and (2) the Exel-Laca algebras  of infinite matrices  with $\{0,1\}$-entries \cite{exe}. It is shown in \cite{fow} that if $E$ is a graph with no sinks and sources, then the $C^*$-algebra $C^*(E)$ is canonically isomorphic to the  Exel-Laca algebra $\mathcal O_{A_E}$, where $A_E$ is the edge matrix of $E$. However, the class of graph $C^*$-algebras and Exel-Laca algebras are differ from each other.

To study both graph $C^*$-algebras and Exel-Laca algebras under one theory,  Tomforde \cite{tom1} introduced the notion of an ultragraph and its associated $C^*$-algebra. Briefly, an ultragraph is a directed graph which allows the range of each edge to be a nonempty set of vertices rather than a singleton vertex. We see in \cite{tom1} that for each binary-valued matrix $A$ there exists an ultragraph $\g_A$ so that the $C^*$-algebra of $\g_A$ is isomorphic to the Exel-Laca algebra of $A$. Furthermore, every graph $C^*$-algebra can be considered as an ultragraph $C^*$-algebra, whereas there is an ultragraph $C^*$-algebra which is not a graph $C^*$-algebra nor an Exel-Laca algebra.

Recently many authors have discussed the algebraic versions of matrix and graph $C^*$-algebras. For example, in \cite{ara}  the algebraic analogue of the Cuntz-Krieger algebra $\mathcal O_A$, denoted by $\mathcal{CK}_A(K)$,  was studied for finite matrix $A$, where  $K$ is a field. Also the Leavitt path algebra $L_K(E)$ of directed graph $E$ was introduced in \cite{abr,abr1} as the algebraic version of graph $C^*$-algebra $C^*(E)$. The class of Leavitt path algebras includes naturally the algebras $\mathcal{CK}_A(K)$ of \cite{ara} as well as the well-known Leavitt algebras $L(1,n)$ of \cite{lea}. More recently, Tomforde defined a new version of Leavitt path algebras with coefficients in a unital commutative ring \cite{tom2}. In the case $R=\mathbb C$, the Leavitt path algebra $L_{\mathbb C}(E)$ is isomorphic to a dense $*$-subalgebra of the graph $C^*$-algebra $C^*(E)$ \cite{tom4}.

The purpose of this paper is to define the algebraic versions of ultragraphs $C^*$-algebras and Exel-Laca algebras. For an ultragraph $\g$ and unital commutative ring $R$, we define the Leavitt path algebra $L_R(\g)$. To study the ideal structure of $L_R(\g)$, we use the notion of quotient ultragraphs from \cite{lar2}. Given an admissible pair $(H,S)$ in $\g$, we define the Leavitt path algebra $L_R(\qg)$ associated to the quotient ultragraph $\qg$ and we prove two kinds of uniqueness theorems, namely the Cuntz-Krieger and the graded-uniqueness theorems, for $L_R(\qg)$. Next we apply these uniqueness  theorems to analyze the ideal structure of $L_R(\g)$. Although the construction of Leavitt path algebra of ultragraph will be similar to that of ordinary graph, we see in Sections 3 and 4 that the analysis of its structure is more complicated. The aim of the definition of ultragraph Leavitt path algebras can be summarized as follows:
\begin{enumerate}
	\item[$\bullet$] Every Leavitt path algebra of a directed graph can be embedded as an subalgebra in a unital ultragraph Leavitt path algebra. Also, the ultragraph Leavitt path algebra $L_{\mathbb C}(\g)$ is isomorphic to a dense $*$-subalgebra of $C^*(\g)$.
	\item[$\bullet$] By using the definition of ultragraph Leavitt path algebras, we give an algebraic version of Exel-Laca algebras.
	\item[$\bullet$] The class of ultragraph Leavitt path algebras is strictly larger than the class of Leavitt path algebras of directed graphs.
\end{enumerate}

The article is organized as follows. We define in Section 2 the Leavitt path algebra $L_R(\g)$ of an ultragraph $\g$ over a unital commutative ring $R$. We continue by considering the definition of quotient ultragraphs of \cite{lar2}. For any admissible pair $(H,S)$ in an ultragraph $\g$, we associate the Leavitt path algebra $L_R\big(\qg\big)$ to the quotient ultragraph $\qg$ and we see that the Leavitt path algebras $L_R(\g)$ and $L_R (\qg\big)$ have a similar behavior in their structure. Next, we prove versions of the graded and Cuntz-Krieger uniqueness theorems for $L_R\big(\qg\big)$ by approximating $L_R\big(\qg\big)$ with $R$-algebras of finite graphs.

By applying the graded-uniqueness theorem in Section 3, we give a complete description of basic graded ideals of $L_R(\g)$ in terms of admissible pairs in $\g$. In Section 4, we use the Cuntz-Krieger uniqueness theorem to show that an ultragraph $\g$ satisfies Condition (K) if and only if every basic ideal in $L_R(\g)$ is graded.

In Section 5, we generalize the algebraic Cuntz-Krieger algebra $\mathcal{CK}_A(K)$ of \cite{ara}, denoted by $\el$, for every infinite matrix $A$ with entries in $\{0,1\}$ and every unital commutative ring $R$. In the case $R=\mathbb C$,  the Exel-Laca $\mathbb C$-algebra $\mathcal{EL}_{A}(\mathbb C)$ is isomorphic to dense $*$-subalgebra of  $\mathcal O_A$. We prove that the class of Leavitt path algebras of ultragraphs contains  the Leavitt path algebras as well as the algebraic Exel-Laca algebras. Furthermore, we give an ultragraph $\g$ such that the Leavitt path algebra $L_R(\g)$  is neither a Leavitt path algebra of graph nor an Exel-Laca $R$-algebra.

\section{Leavitt path algebras}

In this section, we follow the standard constructions of \cite{abr} and \cite{tom1} to define the Leavitt path algebra of an ultragraph. Since the quotient of ultragraph is not an ultragraph, we will have to define the Leavitt path algebras of quotient ultragraphs and prove the uniqueness theorems for them  to characterize the ideal structure in Section 3.

\subsection{Ultragraphs}
Recall from \cite{tom1} that an \emph{ultragraph} $\g=(G^0,\gl,\rg,\sg)$ consists of a set of vertices $G^0$, a set of edges $\gl$, the source map $\sg: \gl \rightarrow G^0$ and the range map $\rg:\gl \rightarrow \mathcal{P}(G^0)\setminus \{\emptyset\}$, where $\mathcal{P}(G^0)$ is the collection of all subsets of $G^0$. Throughout this work, ultragraph $\g$  will be assumed to be countable in the sense that both $G^0$ and $\gl$ are countable.

For a set $X$, a subcollection $\mathcal C$ of $\mathcal P(X)$ is said to be  \emph{lattice} if $\emptyset\in\mathcal C$ and it is closed under the set operations $\cap$ and $\cup$. An \emph{algebra} is a lattice $\mathcal C$ such that $A\setminus B\in \mathcal C$ for all $A,B\in\mathcal C$. If $\g$ is an ultragraph, we write $\go$ for the algebra in $\mathcal P(G^0)$ generated by $\big\{\{v\},\rg(e):v\in G^0, e\in\gl\big\}.$

A \emph{path} in  ultragraph $\g$ is a sequence $\alpha=e_1e_2\cdots e_n$ of edges with $\sg(e_{i+1})\in \rg(e_i)$ for $1\leq i\leq n-1$ and we say that the path $\alpha$ has length $|\alpha|:=n$. We write $\g^n$ for the set of all paths of length $n$ and $\text{Path}(\g):=\bigcup_{n = 0}^\infty  \g^n$ for the set of finite paths. We may extend the maps $\rg$ and $\sg$ on $\text{Path}(\g)$ by setting $\rg(\alpha):=\rg(e_n)$ and $\sg(\alpha):=\sg(e_1)$ for $|\alpha|\ge 1$ and $\rg(A)=\sg(A):=A$ for  $A\in\g^0$.  For every edge $e$, We say that $e^*$ is the \emph{ghost edge} associated to $e$. The  following definition is the algebraic version of \cite[Definition 2.7]{tom1}.
\vspace{0.1in}\\
\noindent
{\bf Definition 2.1.} Let $\g$ be an ultragraph and let $R$ be a unital commutative ring. A \emph{Leavitt $\g$-family} in an $R$-algebra $A$ is a set $\{p_A,s_e,s_{e^*}:A\in\go\,\,\mathrm{and}\,\,e\in\gl\}$ of elements in $A$ such that
\begin{enumerate}
	\item[(1)]$p_\emptyset=0$, $p_A p_B=p_{A\cap B}$ and $p_{A\cup B}=p_A +p_B-p_{A\cap B}$ for all $A,B\in\go$;
	\item[(2)] $p_{\sg(e)}s_e=s_ep_{\rg(e)}=s_e$ and $p_{\rg(e)}s_{e^*}=s_{e^*}p_{\sg(e)}=s_{e^*}$ for all $e\in\gl$;
	\item[(3)] $s_{e^*}s_f=\delta_{e,f}p_{\rg(e)}$ for all $e,f\in\gl$;
	\item[(4)] $p_v=\sum_{\sg(e)=v}s_e s_{e^*}$ for every vertex $v$ with $0<|\sg^{-1}(v)|<\infty$,
\end{enumerate}
where $p_v$ denotes $p_{\{v\}}$. The $R$-algebra generated by the Leavitt $\g$-family $\{s,p\}$ is denoted by $L_R(s,p)$.
\medskip

We say that the Leavitt $\g$-family $\{s,p\}$ is \emph{universal}, if $B$ is an $R$-algebra and $\{S,P\}$ is a Leavitt $\g$-family in $B$, then there exists an algebra homomorphism $\phi:L_R(s,p)\rightarrow  B$ such that $\phi(p_A)=P_A$, $\phi(s_e)=S_e$ and $\phi(s_{e^*})=S_{e^*}$ for every $A\in\go$ and every $e\in\gl$. The \emph{Leavitt path algebra of $\g$ with coefficients in $R$}, denoted by $L_R(\g)$, is the (unique up to isomorphism)  $R$-algebra generated by a universal Leavitt $\g$-family.

\subsection{Quotient ultragraphs}
 We will use the notion of quotient ultragraphs and we generalize the definition of Leavitt path algebras for quotient ultragraphs.
\vspace{0.1in}\\
\noindent
{\bf Definition 2.2.} \cite[Definition 3.1]{tom3}
Let $\g=(G^0,\gl,\rg,\sg)$ be an ultragraph. A subcollection $H\subseteq\go$ is called \emph{hereditary} if satisfying the following conditions:
\begin{enumerate}
	\item $\{\sg(e)\}\in H$ implies $\rg(e)\in H$ for all $e\in\gl$.
	\item $A\cup B\in H$ for all $A,B\in H$.
	\item $A\in H$, $B\in\go$ and $B\subseteq A$, imply $B\in H$.
\end{enumerate}
Also, $H\subseteq\go$ is called \emph{saturated} if for any $v\in G^0$ with $0<|\sg^{-1}(v)|<\infty$, we have
$$\big\{\rg(e):e\in\gl\,\,\mathrm{and}\,\, s_{\g}(e)=v \big\}\subseteq H\,\,\,\,\mathrm{implies}\,\, \,\, \{v\}\in H.$$

For a saturated hereditary subcollection $H\subseteq\go$, the \emph{breaking vertices} of $H$ is denoted by
$$B_H:=\Big\{ v\in G^0: \big|\sg^{-1}(v)\big|=\infty \,\,\mathrm{but}~~ 0<\big|\sg^{-1}(v)\cap \{e:\rg(e)\notin H\}\big|<\infty\Big\}.$$
An \emph{admissible pair in $\g$} is a pair $(H,S)$ of a saturated hereditary set $H\subseteq \go$ and some $S\subseteq B_H$.

In order to define the quotient of ultragraphs we need to recall and introduce some notations from \cite[Section 3]{lar2}. Let $(H,S)$ be an admissible pair in $\g$. For each $A\in\go$, we denote $\overline{A}:=A\cup \{w':w\in A\cap (B_H\setminus S)\}$, where $w'$ denotes another
copy of $w$. Consider the ultragraph $\overline{\g}:=(\overline{G}^0,\overline{\g}^1,\overline{r},\overline{s})$, where $\overline{\g}^1:=\gl$, $\overline{G}^0:=G^0\cup \{w':w\in B_H\setminus S\}$ and the maps $\overline{r}$, $\overline{s}$ are defined by $\overline{r}(e):=\overline{r_{\mathcal G}(e)}$ and
\[\overline{s}(e):=
\left\{
\begin{array}{ll}
s_{\mathcal G}(e)' & \mathrm{if}~s_{\mathcal G}(e)\in  B_H\setminus S ~\mathrm{and}~r_{\mathcal G}(e)\in H,\\
s_{\mathcal G}(e) & \text{otherwise},\\
\end{array}
\right.
\]
for every $e\in\overline{\g}^1$, respectively. We write $\overline{\g}^0$ for the algebra generated by the sets $\{v\}$, $\{w'\}$ and $\overline{r}(e)$.
\vspace{0.1in}\\
\noindent
{\bf Lemma 2.3.} {\cite[Lemma 3.5]{lar2}}
{\it Let $(H,S)$ be an admissible pair in an ultragraph $\g$ and let $\sim$ be a relation on $\overline{\mathcal{G}}^0$ defined by $A\sim B$ if and only if there exists $V\in H$ such that $A\cup V=B\cup V$.	Then $\sim$ is an equivalence relation  on $\overline{\mathcal{G}}^0$ and the operations
$$[A]\cup[B]:=[A\cup B], ~ [A]\cap [B]:=[A\cap B] ~ \mathrm{and}~ [A]\setminus [B]:= [A\setminus B]$$ 
are well-defined on the equivalence classes $\{[A]: A\in \overline{\mathcal{G}}^0\}$.}
\medskip

We usually denote $[v]$ instead of $[\{v\}]$ for every  $v\in \overline{G}^0$. For $A,B\in\overline{\g}^0$, we write $[A]\subseteq [B]$ whenever $[A]\cap [B]=[A]$. The set $\bigcup_{A\in H}A$ is denoted by $\bigcup H$.
\vspace{0.1in}\\
\noindent
{\bf Definition 2.4.} \cite[Definition 3.6]{lar2} Let $(H,S)$ be an admissible pair in  $\g$.  The \emph{quotient ultragraph of $\g$ by $(H,S)$} is the quadruple $\qg:=(\qGo,\qgl,r,s)$, where
$$\qGo:=\big\{[v]:v\in {G}^0\setminus \bigcup H\big\}\cup\big\{[w']:w\in B_H\setminus S\big\},$$
$$\Phi(\mathcal{G}^1):=\big\{e\in{\mathcal{G}}^1:{ r_{\mathcal{G}}}(e)\notin H\big\},$$
and $s: \qgl \rightarrow \qGo$ and $ r:\qgl \rightarrow \{[A]:A\in\overline{\mathcal{G}}^0\}$ are the maps defined by $s(e):=[s_{\mathcal{G}}(e)]$ and $r(e):=[\overline{r_{\mathcal{G}}(e)}]$
for every $e\in\qgl$, respectively.
\vspace{0.1in}\\
\noindent
{\bf Lemma 2.5.}
{\it If $\qg$ is a quotient ultragraph, then
$$\qgo=\bigg\{\bigcup\limits_{j=1}^{k}\bigcap\limits_{i=1}^{n_j}A_{i,j}\setminus B_{i,j}:A_{i,j},B_{i,j}\in\qGo\cup\big\{r(e):e\in  \qgl\big\}\bigg\},$$
where $\qgo$  is the smallest algebra in $\{[A]:A\in\overline{\g}^0\}$ containing $$\Big\{[v],[w']:v\in {G}^0\setminus\bigcup H,w\in B_H\setminus S\Big\}\cup\Big\{r(e): e\in\qgl\Big\}.$$}
\noindent
{\it Proof.} We denote by $X$ the right hand side of the above equality. It is clear that $X\subseteq\qgo$, because $\qgo$ is an algebra generated by the elements $[v]$,$[w']$ and $r(e)$. For the reverse inclusion, we note that $X$ is a lattice. Furthermore, one can show that $X$ is closed under the operation $\setminus$. Thus
$X$ is an algebra contains $[v]$,$[w']$ and $r(e)$ and consequently $\qgo\subseteq X$.
\hfill$\Box$
\medskip

\noindent
{\bf Remark 2.6.} If $A,B\in\go$, then $\overline{A\cup B}=\overline{A}\cup \overline{B},\, \overline{A\cap B}=\overline{A}\cap \overline{B}$ and $\overline{A\setminus B}=\overline{A}\setminus\overline{B}$. Thus, by applying an analogous lemma of Lemma 2.5 for $\go$ and $\overline{\g}^0$, we deduce that $\overline{A}\in\overline{\g}^0$ for all $A\in\go$.
One can see that
\begin{multline*}
\overline{\g}^0=\big\{\overline A\cup F_1\cup F_2:A\in\go,~ F_1~\text{and}~F_2~\text{are}~\text{finite~subsets~of}~G^0\\
\text{and}~\{w':w\in B_H\},~\text{respectively}\big\}.
\end{multline*}
For example we have
$$\overline{A}\setminus\{v\}=\overline{A\setminus\{v\}}\cup\Big(\overline{\{v\}}\setminus\{v\}\cap\overline{A}\setminus A\Big),$$
and
$$\overline{A}\setminus\{w'\}=\overline{A\setminus\{w\}}\cup\Big(A\cap\{w\}\Big).$$
Furthermore, it follows from Lemma 2.3 that $\qgo=\big\{[A]:A\in\overline{\g}^0\big\}.$
\medskip

\noindent
{\bf Remark 2.7.}
The hereditary property of $H$ and Remark 2.6 imply that $A\sim B$ if and only if both $A\setminus B$ and $B\setminus A$ belong to $H$.
\medskip

Similar to ultragraphs, a \emph{path} in $\qg$ is a sequence $\alpha:=e_1e_2\cdots e_n$ of edges in $\qgl$ such that $s(e_{i+1})\subseteq r(e_i)$ for $1\leq i\leq n-1$. We say the path $\alpha$ has length $|\alpha|:=n$ and we consider the elements in $\qgo$ to be paths of length zero. We denote by $\mathrm{Path}(\qg)$, the union of paths with finite length. The maps $r$ and $s$ can be naturally extended on $\mathrm{Path}(\qg)$. Let $\qgl^*$ be the set of \emph{ghost edges} $\{e^*:e\in\qgl \}$. We also define the \emph{ghost path} $\alpha^*:=e_n^*e_{n-1}^*\cdots e_1^*$ for every $\alpha=e_1e_2\cdots e_n\in\mathrm{Path}(\qg)\setminus\qgo$ and $[A]^*:=[A]$ for every $[A]\in\qgo$.

Using Theorem 3.4, we define the Leavitt path algebra of a quotient ultragraph $\qg$ which is corresponding to the quotient $R$-algebra $L_R(\g)/I_{(H,S)}$. We use this concept to characterize the ideal structure of  $L_R(\g)$ in Section 3.  The  following definition is the algebraic version of \cite[Definition 3.8]{lar2}.
\vspace{0.1in}\\
\noindent
{\bf Definition 2.8.} Let $\qg$ be a quotient ultragraph and let $R$ be a unital commutative ring. A \emph{Leavitt $\qg$-family} in an $R$-algebra $A$ is a set $\{q_{[A]},t_e,t_{e^*}:[A]\in\qgo\,\,\mathrm{and}\,\,e\in\qgl\}$ of elements in $A$ such that
\begin{itemize}
	\item[(1)]$q_{[\emptyset]}=0$, $q_{[A]} q_{[B]}=q_{[A]\cap [B]}$ and $q_{[A]\cup [B]}=q_{[A]} +q_{[B]}-q_{[A]\cap [B]}$;
	\item[(2)] $q_{s(e)}t_e=t_eq_{r(e)}=t_e$ and $q_{r(e)}t_{e^*}=t_{e^*}q_{s(e)}=t_{e^*}$;
	\item[(3)] $t_{e^*}t_f=\delta_{e,f}q_{r(e)}$;
	\item[(4)] $q_{[v]}=\sum_{s(e)=[v]}t_e t_{e^*}$ for every $[v]\in\qGo$ with $0<|s^{-1}([v])|<\infty$.
\end{itemize}
The $R$-algebra generated by the Leavitt $\qg$-family $\{t,q\}$ is denoted by $L_R(t,q)$. The \emph{Leavitt path algebra of $\qg$ with coefficients in $R$}, denoted by $\lqg$, is the (unique up to isomorphism)  $R$-algebra generated by a universal Leavitt $\qg$-family (the definition of universal Leavitt $\qg$-family is similar to ultragraph case).
\medskip

If $(H,S)=(\emptyset,\emptyset)$, then we have $[A]=\{A\}$ for each $A\in \qgo$. In this case, every Leavitt $\g/(\emptyset,\emptyset)$-family is a Leavitt $\g$-family and vice versa. So, we can consider the ultragraph Leavitt path algebra $\lrg$ as $L_R(\g/(\emptyset,\emptyset))$.

Let $R$ be a unital commutative ring. For a nonempty set $X$, we write $w(X)$ for the set of words $w:=w_1w_2\cdots w_n$ from the alphabet $X$. The \emph{free $R$-algebra} generated by $X$ is denoted by $\mathbb F_R(w(X))$. For definition of the free $R$-algebra we refer the reader to \cite[Section 2.3]{pask}.

Now, we show that for every quotient ultragraph $\qg$, there exists a universal Leavitt $\qg$-family. Suppose that  $X:=\qgo\cup\qgl\cup\qgl^*$ and $I$ is the ideal of the free $R$-algebra $\mathbb F_R(w(X))$ generated by the union of the following sets:
\begin{itemize}
	\item $\big\{[\emptyset],\, [A][B]-[A\cap B],\, [A\cup B]+[A\cap B]-([A] +[B]):[A],[B]\in\qgo\big\},$
	\item $\big\{e-s(e)e,e-er(e),e^*-e^*s(e),e^*-r(e)e^*:e\in\qgl\big\},$
	\item $\big\{e^*f-\delta_{e,f}r(e):e,f\in\qgl\big\},$
	\item $\big\{v-\sum_{s(e)=[v]}ee^*:0<|s^{-1}([v])|<\infty\big\}$.
\end{itemize}
If $\pi:\mathbb F_R(w(X))\rightarrow \mathbb F_R(w(X))/I$ is the quotient map, then it is easy to check that the collection $\{\pi([A]),\pi(e),\pi(e^*):[A]\in\qgo,e\in\qgl\}$ is a Leavitt $\qg$-family. For our convenience, we denote $q_{[A]}:=\pi(A)$, $t_e:=\pi(e)$ and $t_{e^*}:=\pi(e^*)$ for every $[A]\in\qgo$ and $e\in\qgl$, and we show that the Leavitt $\g$-family $\{t,q\}$ has the universal property.

Assume that $\{T,Q\}$ is a Leavitt $\qg$-family in an $R$-algebra $A$. If we define $f:X\rightarrow A$ by $f([A])=P_{[A]}$, $f(e)=T_e$ and $f(e^*)=T_{e^*}$, then by  \cite[Proposition 2.6]{pask}, there is an $R$-algebra homomorphism $\phi_f:\mathbb F_R(w(X))\rightarrow A$ such that $\phi_f\big|_X=f$. Since $\{T,Q\}$ is a Leavitt $\qg$-family, $\phi_f$ vanishes on the generator of $I$. Hence, we can define an $R$-algebra homomorphism $\phi:\mathbb F_R(w(X))/I\rightarrow A$ by $\phi\big(\pi(x)\big)=\phi_f(x)$ for every $x\in\mathbb F_R(w(X))$. Furthermore, we have $\phi(q_{[A]})=Q_{[A]}$, $\phi(t_e)=T_e$ and $\phi(t_{e^*})=T_{e^*}$.

From now on we denote the universal Leavitt $\g$-family and $\qg$-family by
$\{s,p\}$ and $\{t,q\}$, respectively. So, we suppose that $\{s,p\}$ and $\{t,q\}$ are the canonical generators of $\lrg$ and $\lqg$, respectively.
\vspace{0.1in}\\
\noindent
{\bf Theorem 2.9.}
{\it Let $\qg$ be a quotient ultragraph and let $R$ be a unital commutative ring. Then $\lqg$ is of the form
	$$\mathrm{span}_R\big\{t_\alpha q_{[A]}t_{\beta^*}\in\lqg: \alpha,\beta\in \mathrm{Path}\big(\qg\big), r(\alpha)\cap [A]\cap r(\beta)\ne[\emptyset]  \big\}.$$
	Furthermore, $\lqg$ is a $\mathbb{Z}$-graded ring by the grading
	$$\lqg_n=\mathrm{span}_R\big\{t_{\alpha}q_{[A]}t_{\beta^*}\in\lqg:|\alpha|-|\beta|=n\big\} \hspace{.5cm} (n\in\mathbb{Z}).$$}
\noindent
{\it Proof.}
If $t_\alpha q_{[A]}t_{\beta^*}, t_\mu q_{[B]}t_{\nu^*}\in\lqg$, then
\[(t_\alpha q_{[A]}t_{\beta^*})(t_\mu q_{[B]}t_{\nu^*})=
\left\{
\begin{array}{ll}
t_{\alpha\mu^{'}}q_{[B]}t_{\nu^*} & \mathrm{if}\,\, \mu=\beta\mu^{'}\,\,\mathrm{and}\,\, s(\mu^{'})\subseteq A\cap r(\alpha), \\
t_\alpha q_{[A]\cap r(\beta)\cap [B]}t_{\nu^*} & \mathrm{if}\,\,\mu=\beta,\\
t_\alpha q_{[A]}t_{(\nu\beta^{'})^*} & \mathrm{if}\,\, \beta=\mu\beta^{'}\,\,\mathrm{and}\,\, s(\beta^{'})\subseteq B\cap r(\nu),\\
0 & \mathrm{otherwise},
\end{array}
\right.
\]
which proves the first statement. Let $X:=\qgo\cup\{e,e^*:e\in\qgl\}$ and consider $\lqg=\mathbb F_R(w(X))/I$ as defined before. If we define a degree map $d:X\rightarrow\mathbb Z$ by $d([A])=0$, $d(e)=1$ and $d(e^*)=-1$ for every $[A]\in\qgo$ and $e\in\qgl$, then by \cite[Proposition 2.7]{pask}, $\lqg$ is a $\mathbb Z$-graded ring with the grading 
$$\lqg_n=\mathrm{span}_R\big\{t_{\alpha}q_{[A]}t_{\beta^*}\in\lqg:|\alpha|-|\beta|=n\big\}.$$ 
\hfill$\Box$
\medskip

\noindent
{\bf Theorem 2.10.}
{\it Let $\g$ be an ultragraph and let $R$ be a unital commutative ring. Then for all nonempty sets $A\in\go$ and all $r\in R\setminus\{0\}$, the elements $rp_A$ of $\lrg$ are nonzero. In particular, $rs_e\ne0$ and  $rs_{e^*}\ne0$ for every $A\in\go\setminus\{\emptyset\}$ and every $r\in R\setminus\{0\}$.}
\vspace{0.05in}\\
\noindent
{\it Proof.} By the universality, it suffices to generate a Leavitt $\g$-family $\{\tilde s,\tilde p\}$ in which  $r\tilde p_A\ne0$ for every nonempty set $A\in\go$ and every $r\in R\setminus\{0\}$. For each $e\in\gl$, we define a disjoint copy $Z_e:=\bigoplus R$, where the  direct sum is taken over countably many copies of $R$ and for each $v\in G^0$, let
\[Z_v:=
\left\{
\begin{array}{ll}
\bigoplus\limits_{s(e)=v}Z_e & \mathrm{if}\,\, |s^{-1}_{\g}(v)|\ne0,\\ \\
~~\bigoplus R & \mathrm{if}\,\, |s^{-1}_{\g}(v)|=0,
\end{array}
\right.
\]
where $\bigoplus R$ is another disjoint copy for each $v$. For every $\emptyset\ne A\in\go$, define $P_A:\bigoplus_{v\in A}Z_v\rightarrow \bigoplus_{v\in A}Z_v$ to be the identity map. Also, for each $e\in\gl$ choose an isomorphism $S_e:\bigoplus_{v\in r_{\g}(e)}Z_v\rightarrow Z_e$ and let $S_{e^*}:=S_e^{-1}:Z_e\rightarrow \bigoplus_{v\in r_{\g}(e)}Z_v$. Now if $Z:=\bigoplus_{v\in G^0}Z_v$, then we naturally extend all $P_A, S_e,$ $S_{e^*}$ to homomorphisms $\tilde p_A, \tilde s_e,\tilde s_{e^*}\in \mathrm{Hom}_R(Z,Z)$, respectively by setting the yet undefined components to zero. It is straightforward to verify that $\{\tilde s,\tilde p\}$ is a Leavitt $\g$-family in $\mathrm{Hom}_{R}(Z,Z)$ such that $r\tilde p_A\ne0$ for every $\emptyset\ne A\in\go$ and every $r\in R\setminus\{0\}$.
\hfill$\Box$
\medskip

Note that we cannot follow the argument of Theorem 2.10 to show that $rq_{[A]}\ne0$. For example, suppose that $\g$ is the  ultragraph
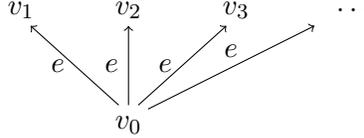
\begin{center}
	\begin{tikzpicture} [x=37pt,y=30pt]
	\draw [->] (-.1,.25) -- node[ left] {$e$} (-1,1.25);
	\draw (-1.1,1.45) node(x)  {$v_1$} (-1.1,1.45);
	\draw [->] (0,.25) -- node[ left] {$e$} (0,1.25);
	\draw (0,1.45) node(x)  {$v_2$} (0,1.45);
	\draw [->] (.1,.25) -- node[ left] {$e$} (1,1.25);
	\draw (1.1,1.45) node(x)  {$v_3$} (1.1,1.45);
	\draw [->] (.2,.20) -- node[above ] {$e$} (1.9,1.25);
	\draw (2.3,1.45) node(x)  {$\ldots$} (2.3,1.45);
	\draw (0,0) node(x)  {$v_0$} (0,0);
	\end{tikzpicture}
\end{center}
and let $H$ be the collection of all finite subsets of $\{v_1,v_2,\ldots\}$, which is a hereditary and saturated subcollection of $\go$. If we consider the quotient ultragraph $\g/(H,\emptyset)$, then $\{[\emptyset]\ne[v]:[v]\subseteq r(e)\}=\emptyset$. So  we can not define the idempotent $q_{r(e)}:\bigoplus_{[v]\subseteq r(e)}Z_{[v]}\rightarrow \bigoplus_{[v]\subseteq r(e)}Z_{[v]}$ as in the proof of Theorem 2.10. In Section 3 we will solve this problem.
\vspace{0.05in}\\
\noindent
{\bf Remark 2.11.} Every directed graph $E=(E^0,E^1,r,s)$ can be considered as an ultragraph $\g=(G^0,\gl,\overline{r}_{\g},s_{\g})$, where $G^0:=E^0$, $\gl:=E^1$ and the map $\overline r_{\g}:\gl\rightarrow\mathcal{P}(G^0)\setminus \{\emptyset\}$ is defined by $\overline r_{\g}(e)=\{r(e)\}$ for every $e\in\gl$. In this case, the algebra $\go$ is the collection of all finite subsets of $G^0$. The Leavitt path algebra $L_R(E)$ is naturally isomorphic to $L_R(\g)$ (see \cite{abr,tom2} for more details about the Leavitt path algebras of directed graphs). So the class of ultragraph Leavitt path algebras contains the class of  Leavitt path algebras of directed graphs.
\medskip

\noindent
{\bf Lemma 2.12.}
{\it Let $\g$ be an ultragraph and let $R$ be a unital commutative ring. Then $L_R(\g)$ is unital if and only if $G^0\in\go$ and in this case $1_{L_R(\g)}=p_{G^0}$.}
\vspace{0.05in}\\
\noindent
{\it Proof.} If $G^0\in\go$, then the Relations of Definition 2.1 imply that $p_{G^0}$ is a unit for $L_R(\g)$.

Conversely, suppose that $L_R(\g)$ is unital and write
$$1_{L_R(\g)}=\sum_{k=1}^n r_ks_{\alpha_k}p_{A_k}s_{\beta^*_k},$$
where $r_k\in R, A_k\in\go$ and $\alpha_k,\beta_k\in\mathrm{Path}(\g)$. Let $A:=\bigcup_{k=1}^ns(\alpha_k)\in\go$. If $G^0\notin\go$, then we can choose an element $v\in G^0\setminus A$ and derive
$$p_v=p_v\cdot1_{L_R(\g)}=\sum_{k=1}^nr_kp_vs_{\alpha_k}p_{A_k}s_{\beta^*_k}= \sum_{k=1}^nr_kp_v p_{s_{\g}(\alpha_k)}s_{\alpha_k}p_{A_k}s_{\beta^*_k}=0,$$
this  contradicts  Theorem 2.10 and it follows $G^0\in\go$, as desired.
\hfill$\Box$
\medskip

We note that, for directed graph $E$, the algebra $L_R(E)$ is unital if and only if $E^0$ is finite. If $E^0$ is infinite, then we can define an ultragraph $\g$ associated to $E$ such that $L_R(\g)$ is unital and $L_R(E)$ is an embedding subalgebra of $L_R(\g)$. More precisely, Consider the ultragraph $\g=(G^0,\gl,r_{\g},s_{\g})$ where $G^0=E^0\cup\{v\}$, $\gl=E^1\cup\{e\}$, $s_{\g}(e)=v$ and $r_{\g}(e)=E^0$. Since $r_{\g}(e)\cup\{v\}=G^0\in\go$, by Lemma 2.12, $L_R(\g)$ is unital. Define $S_e=s_e$ and $P_v=p_v$ for $e\in E^1$ and $v\in E^0$, respectively. It is straightforward to see that $\{S,P\}$ is a Leavitt $E$-family for directed graph $E$. The result now follows by \cite[Theorem 5.3]{tom2}.

\subsection{Uniqueness Theorems}

Let $\qg$ be a quotient ultragraph. We prove the graded and Cuntz-Krieger uniqueness theorems for $\lqg$ and $L_R(\g)$. We do this  by approximating the Leavitt path algebras of quotient ultragraphs  with the Leavitt path algebras of  finite graphs.
Our proof in this section is standard (see \cite[Section 4]{lar2}), and we give the details for simplicity of further results of the paper.

A vertex $[v]\in \qGo$ is called a \emph{sink} if $s^{-1}([v])=\emptyset$ and $[v]$ is called an \emph{infinite emitter} if $|s^{-1}([v])|=\infty$. A \emph{singular vertex} is a vertex that is either a sink or an infinite emitter. The set of all singular vertices is denoted by $\qsg$.

Let $F\subseteq \qsg \cup \qgl$ be a finite subset and write $F^0:=F\cap \qsg$ and $F^1:=F\cap \qgl=\{e_1,\ldots,e_n\}$. Following \cite{lar2}, we construct a finite graph $E_F$ as follows. First, for every $\omega=(\omega_1,\ldots, \omega_n)\in \{0,1\}^n\setminus \{0^n\}$, we define  $r(\omega):=\bigcap_{\omega_i=1}r(e_i)\setminus \bigcup_{\omega_j=0}r(e_j)$ and $R(\omega):=r(\omega)\setminus \bigcup_{[v]\in F^0}[v]$ which belong to $\qgo$. We also set
\begin{multline*}
\Gamma_0:=\{\omega\in \{0,1\}^n\setminus\{0^n\}: $ there are  vertices $ [v_1],\ldots,[v_m]\\
$ such  that $ R(\omega)=\bigcup_{i=1}^m[v_i] $ and $ \emptyset\neq s^{-1}([v_i])\subseteq F^1 $ for $ 1\leq i\leq m\},
\end{multline*}
and
\[\Gamma_F:=\left\{\omega\in \{0,1\}^n\setminus \{0^n\}: R(\omega)\neq [\emptyset] \mathrm{~and ~} \omega\notin \Gamma_0 \right\}.\]
Now we define the finite graph $E_F=(E_F^0,E_F^1,r_F,s_F)$, where
\begin{align*}
E_F^0:=&F^0 \cup F^1 \cup \Gamma_F,\\
E_F^1:=&\left\{(e,f)\in F^1\times F^1: s(f)\subseteq r(e) \right\}\\
&\cup \left\{(e,[v])\in F^1\times F^0: [v]\subseteq r(e) \right\}\\
&\cup \left\{(e,\omega)\in F^1\times \Gamma_F: \omega_i=1 \mathrm{~~whenever~~} e=e_i \right\},
\end{align*}
with $s_F(e,f)=s_F(e,[v])=s_F(e,\omega)=e$ and $r_F(e,f)=f$, $r_F(e,[v])=[v]$, $r_F(e,\omega)=\omega$.
\vspace{0.1in}\\
\noindent
{\bf Lemma 2.13.}
{\it Let $\qg$ be a quotient ultragraph and let $R$ be a unital commutative ring. Then we have the following assertion:\\
	(i) For every finite set $F\subseteq \qsg \cup \qgl$, the elements
	$$\begin{array}{lll}
	P_e:=t_et_{e^*}, &P_{[v]}:=q_{[v]}(1-\sum\limits_{e\in F^1}t_et_{e^*}), &P_\omega:=q_{R(\omega)}(1-\sum\limits_{e\in F^1}t_et_{e^*}), \\
	S_{(e,f)}:=t_eP_f, &S_{(e,[v])}:=t_e P_{[v]}, &S_{(e,\omega)}:=t_e P_\omega,\\
	S_{(e,f)^*}:=P_ft_{e^*}, &S_{(e,[v])^*}:=P_{[v]}t_{e^*}, &S_{(e,\omega)^*}:=P_\omega t_{e^*},
	\end{array}$$
	form a Leavitt $E_F$-family which generates the subalgebra of $\lqg$ generated by $\left\{q_{[v]},t_e,t_{e^*}:[v]\in F^0,e\in F^1\right\}$.\\
	(ii) For $r\in R\setminus\{0\}$, if $rq_{[A]}\ne 0$ for all $[A]\ne[\emptyset]$ in $\qgo$, then  $rP_z\ne0$ for all $z\in E^0_F$. In this case, we have 
	$$L_R(E_F)\cong L_R(S,P)=\mathrm{Alg}\{q_{[v]},t_e,t_{e^*}\in\lqg:[v]\in F^0,e\in F^1\}.$$}
\noindent
{\it Proof.} The statement (i) follows from the fact that $\{t,q\}$ is a Leavitt
$\qg$-family, or see the similar \cite[Proposition 4.2]{lar2}.

For (ii), fix $r\in R\setminus \{0\}$. If $rq_{[A]}\ne0$ for every $[A]\in\qgo\setminus\{[\emptyset]\}$, then $rt_e\ne 0$ and $rt_{e^*}\ne0$ for every edge $e$. Thus  $rP_e\ne0$ for every edge $e\in F^1$. Let $[v]\in F^0$. If $[v]$ is a sink, then $rP_{[v]}=rq_{[v]}\ne0$. If $[v]$ is an infinite emitter, then there is $f\in\qgo\setminus F^1$ such that $s(f)=[v]$. In this case, we have $rP_{[v]}t_f=rq_{[v]}t_f=t_f\ne0$. Therefore  $rP_{[v]}\ne0$ for all $[v]\in F^0$. Moreover, for each $\omega\in\Gamma_F$, there is a vertex $[v]\subseteq R(\omega)$ such that either $[v]$ is a sink or there is an edge $f\in\qgl\setminus F^1$ with $s(f)=[v]$. In the former case, we have $q_{[v]}(rP_\omega)=rq_{[v]}\ne0$ and in the later, $t_{f^*}(rP_\omega)=rt_{f^*}\ne0$. Thus all $rP_\omega$ are nonzero. Consequently, $rP_z\ne0$ for every vertex $z\in E^0_F$.

Now we show that $L_R(E_F)\cong L_R(S,P)$. Note that for $z\in E_F^0$ and $g\in E_F^1$, the degree of $P_z$, $S_g$ and $S_{g^*}$ as elements in $\lqg$ are 0, 1 and -1, respectively. So $L_R(S,P)$ is a graded subalgebra of $\lqg$ with the grading 
$$L_R(S,P)_n:=L_R(S,P)\cap L_R(\qg)_n.$$
Let $\{\tilde s,\tilde p\}$ be the canonical generators of $L_R(E_F)$. By the universal property, there is an $R$-algebra homomorphism $\pi:L_R(E_F)\rightarrow L_R(S,P)$ such that $\pi(r\tilde p_z)=rP_z\ne0$, $\pi(\tilde s_g)=S_g$ and $\pi(\tilde s_{g^*})=S_{g^*}$ for $z\in E^0_F$, $g\in E^1_F$ and $r\in R\setminus \{0\}$. Since $\pi$ preserves the degree of generators, the graded-uniqueness theorem for graphs \cite[Theorem 5.3]{tom2} implies that $\pi$ is injective. As $\pi$ is also surjective, we conclude that $L_R(E_F)$ is isomorphic to $L_R(S,P)$ as $R$-algebras.
\hfill$\Box$
\medskip

\noindent
{\bf Theorem 2.14} (The graded-uniqueness theorem).
{\it Let $\qg$ be a quotient ultragraph and let $R$ be a unital commutative ring. If $T$ is a $\mathbb Z$-graded ring and $\pi:\lqg\rightarrow T$ is  a graded ring  homomorphism with $\pi(rq_{[A]})\ne0$ for all $[A]\ne[\emptyset]$ in $\qgo$ and $r\in R\setminus \{0\}$, then $\pi$ is injective.}
\vspace{0.05in}\\
\noindent
{\it Proof.} Let $\{F_n\}$ be an increasing sequence of finite subsets of $\qsg \cup \qgl$ such that $\cup_{n=1}^\infty F_n=\qsg \cup \qgl$. For each $n$, the graded subalgebra of $\lqg$ generated by $\{q_{[v]},t_e,t_{e^*}:[v]\in F_n^0\,\,\mathrm{and}\,\, e\in F_n^1\}$ is denoted by $X_n$. Since $\pi(rq_{[A]})\ne0$ for all $[A]\in\qgo\setminus\{[\emptyset]\}$ and $r\in R\setminus \{0\}$, by Lemma 2.13, there is an graded isomorphism $\pi_n:L_R(E_{F_n})\rightarrow X_n$. Thus $\pi\circ\pi_n:L_R(E_{F_n})\rightarrow T$ is a graded homomorphism. 

Fix $r\in R\setminus \{0\}$. We show that $\pi\circ\pi_n(r\tilde p_z)=\pi(rP_z)\ne0$ for all $z\in E_{F_n}^0$. Since $\pi(rq_{[A]})\ne0$ for all $[A]\in\qgo\setminus\{[\emptyset]\}$, the elements $\pi(rt_e)$ and $\pi(rt_{e^*})$ are nonzero for every edge $e$. So $\pi(rP_et_e)=\pi(rt_e)\ne0$ and thus $\pi(rP_e)\ne0$. Let $[v]\in F^0$. If $[v]$ is a sink, then $\pi(rP_{[v]})=\pi(rq_{[v]})\ne0$. If $[v]$ is an infinite emitter, then there is $f\in\qgo\setminus F^1$ such that $s(f)=[v]$. In this case, we have $\pi(rP_{[v]}t_f)=\pi(rq_{[v]}t_f)=\pi(t_f)\ne0$, hence $\pi(rP_{[v]})\ne0$. For every $\omega\in\Gamma_F$, there is a vertex $[v]\subseteq R(\omega)$ such that either $[v]$ is a sink or there is an edge $f\in\qgl\setminus F^1$ with $s(f)=[v]$. In the former case, we have $\pi(q_{[v]}(rP_\omega))=\pi(rq_{[v]})\ne0$ and in the later, $\pi(t_{f^*}(rP_\omega))=\pi(rt_{f^*})\ne0$. Thus $\pi\circ\pi_n(r\tilde p_z)\ne0$ for every $z\in E_{F_n}^0$ and $r\in R\setminus\{0\}$. Hence, we may apply the graded-uniqueness theorem for graphs \cite[Theorem 5.3]{tom2} to obtain the injectivity of $\pi\circ\pi_n$.

If $[v]$ is a non-singular vertex, then we have $q_{[v]}=\sum_{s(e)=[v]}t_et_{e^*}$. Furthermore, $q_{[A]\setminus[B]}=q_{[A]}-q_{[A]}q_{[B]}$ for every $[A],[B]\in\qgo$.
Thus, by Lemma 2.5, $L_R(\qg)$ is an $R$-algebra generated by $$\left\{q_{[v]},t_e,t_{e^*}:[v]\in\qsg\,\,\mathrm{and}\,\, e\in\qgl\right\},$$
and so $\cup_{n=1}^\infty X_n=\lqg$. It follows that $\pi$ is injective on $\lqg$, as desired.
\hfill$\Box$
\medskip

\noindent
{\bf Corollary 2.15.}
 Let $\g$ be an ultragraph, $R$ a unital commutative ring and $T$ a $\mathbb Z$-graded ring. If $\pi:L_R(\g)\rightarrow T$ is  a graded ring  homomorphism such that $\pi(rp_A)\ne0$ for all $A\in\go\setminus\{\emptyset\}$ and $r\in R\setminus \{0\}$, then $\pi$ is injective.
\medskip

\noindent
{\bf Definition 2.16.} A \emph{loop} in a quotient ultragraph $\qg$ is a path $\alpha$ with $|\alpha|\geq1$ and $s(\alpha)\subseteq r(\alpha)$. An \emph{exit} for a loop $\alpha_1\cdots\alpha_n$ is an edge $f\in\qgl$ with the property that $s(f)\subseteq r(\alpha_i)$ but $f\ne\alpha_{i+1}$ for some $1\leq i\leq n$, where $\alpha_{n+1}:=\alpha_1$. We say that $\qg$ satisfies \emph{Condition (L)} if every loop $\alpha:=\alpha_1\cdots\alpha_n$ in $\qg$ has  an exit, or $r(\alpha_i)\neq s(\alpha_{i+1})$ for some $1\leq i\leq n$.
\medskip

\noindent
{\bf Theorem 2.17} (The Cuntz-Krieger uniqueness theorem).
{\it Let $\qg$ be a quotient ultragraph satisfying Condition (L) and let $R$ be a unital commutative ring. If $T$ is a ring and $\pi:\lqg\rightarrow T$ is a ring homomorphism such that $\pi(rq_{[A]})\ne0$ for every $[A]\in\qgo\setminus\{[\emptyset]\}$ and $r\in R\setminus \{0\}$, then $\pi$ is injective.}
\vspace{0.05in}\\
\noindent
{\it Proof.} Choose an increasing sequence $\{F_n\}$ of finite subsets of $\qsg \cup \qgl$ such that $\cup_{n=1}^\infty F_n=\qsg \cup \qgl$. Let $X_n$ be the subalgebras of $L_R(\qg)$ as in Theorem 2.14. Since $\pi(rq_{[A]})\ne0$ for all $[A]\in\qgo\setminus\{[\emptyset]\}$ and $r\in R\setminus \{0\}$, by Lemma 2.13, there exists an isomorphism $\pi_n:L_R(E_{F_n})\rightarrow X_n$ for each $n\in\mathbb N$. Furthermore, $\pi\circ\pi_n(r\tilde p_z)\ne0$ for every $z\in E_{F_n}^0$ and $r\in R\setminus\{0\}$. Since $\qg$ satisfies Condition (L), By \cite[Lemma 4.8]{lar2}, all finite graphs $E_{F_n}$ satisfy Condition (L). So, the Cuntz-Krieger uniqueness theorem for graphs \cite[Theorem 6.5]{tom2} implies that $\pi\circ\pi_n$ is injective for $n\geq 1$. Now by the fact $\cup_{n=1}^\infty X_n=L_R(\qg)$, we conclude that $\pi$ is an injective homomorphism.
\hfill$\Box$
\medskip

\noindent
{\bf Corollary 2.18.}
Let $\g$ be an ultragraph satisfying Condition (L), $R$ a unital commutative ring and  $T$  a ring. If $\pi:L_R(\g)\rightarrow T$ is  a  ring  homomorphism such that $\pi(rp_A)\ne0$ for all $A\in\go\setminus\{\emptyset\}$ and $r\in R\setminus \{0\}$, then $\pi$ is injective.

%%%%%%%%%%%%%%%%%%%%%%%%%%%%%%%%%%%%%%%%%%%%%%%%
\section{Basic Graded Ideals}\label{s3}

In this section, we apply the graded-uniqueness theorem for quotient ultragraphs to investigate the ideal structure of $L_R(\g)$. We would like to consider the ideals of $L_R(\g)$ that are reflected in the structure of the ultragraph $\g$. For this, we give the following definition of basic ideals.

Let $(H,S)$ be an admissible pair in an ultragraph $\g$. For any $w\in B_H$, we define the gap idempotent
$$p_w^H:=p_w-\sum_{s(e)=w,~ r(e)\notin H}s_e s_{e^*}.$$

Let $I$ be an ideal in $L_R(\g)$. We write $H_{I}:=\{A\in\go:p_A\in I\}$, which is a saturated hereditary subcollection of $\go$. Also, we set $S_I:=\{w\in B_{H_I}:p^{H_I}_w\in I\}$. We say that the ideal
$I$ is \emph{basic} if the following conditions hold:
\begin{enumerate}
	\item $rp_A\in I$ implies $p_A\in I$ for $A\in\go$ and $r\in R\setminus\{0\}$.
	\item $rp_w^{H_I}\in I$ implies $p_w^{H_I}\in I$ for $w\in B_{H_I}$ and $r\in R\setminus\{0\}$.
\end{enumerate}

For an admissible pair  $(H,S)$ in  $\g$,
the (two-sided) ideal of $L_R(\g)$ generated by the idempotents $\{p_A:A\in H\} \cup \left\{p_w^H:w\in S\right\}$  is denoted by $\I$.
\vspace{0.1in}\\
\noindent
{\bf Lemma 3.1.}(cf. \cite[Lemma 3.9]{lar1})
{\it If $(H,S)$ is an admissible pair in ultragraph $\g$, then
	$$\I:=\mathrm{span}_R\big\{s_\alpha p_As_{\beta^*},s_\mu p_w^H s_{\nu^*}\in\lrg:A\in H\,\,\mathrm{and}\,\,w\in S\big\}$$
	and $\I$ is a graded basic ideal of $L_R(\g)$.}
\vspace{0.05in}\\
\noindent
{\it Proof.} We denote the right-hand side of the above equality  by $J$. The hereditary property of $H$ implies that $J$ is an ideal of $L_R(\g)$ being contained
in $\I$. On the other hand, all generators of $I_{(H,S)}$ belong to $J$ and so we have $I_{(H,S)}=J$. Note that the elements $s_\alpha p_As_{\beta^*}$ and $s_\mu p_w^H s_{\nu^*}$ are homogeneous elements of degrees $|\alpha|-|\beta|$ and $|\mu|-|\nu|$, respectively. Thus $\I$ is a graded ideal.

To show that $I_{(H,S)}$ is a basic ideal, suppose $rp_A\in I_{(H,S)}$ for some $A\in\go$ and $r\in R\setminus \{0\}$ and write
$$rp_A=x:=\sum\limits_{i=1}^{n}r_is_{\alpha_i}p_{A_i}s_{\beta^*_i} +\sum\limits_{j=1}^{m}s_js_{\mu_j}p^H_{w_j}s_{\nu^*_j},$$
where $A_i\in H$, $w_j\in S$ and $r_i,s_j\in R$ for all $i,j$. We first show assertion (1) of the definition of basic ideal in several steps.

\textbf{Step I:} If $A=\{v\}\notin H$ and $v\notin S$, then $0<|s^{-1}_{\g}(v)|<\infty$.

Note that $p_vx=x$, so  the assumption $v\notin \bigcup H\cup S$ and the hereditarity of $H$ imply that we may assume that $|\alpha_i|,|\mu_j|\ge1$ and $s_{\g}(\alpha_i)=s_{\g}(\mu_j)=v$ for every $i,j$ (because a sum like $\sum r_iA_i$ could be zero). Hence $v$ is not a sink. Set $\alpha_i=\alpha_{i,1}\alpha_{i,2}\cdots\alpha_{i,|\alpha|}$. If $|s^{-1}_{\g}(v)|=\infty$, then there is an edge $e\ne\alpha_{1,1},\dots,\alpha_{n,1},\mu_{1,1},\dots,\mu_{m,1}$ with $s_{\g}(e)=v$. So $rs_{e^*}=s_{e^*}(rp_v)=s_{e^*}x=0 $, contradicting Theorem 2.3. Thus $0<|s^{-1}_{\g}(v)|<\infty$.

\textbf{Step II:} If $ A\notin H$, then there exists $v\in A$ such that $\{v\}\notin H$.

Let $rp_A=x$. Assume first that $|\mu_j|=0$ for some $j$. Then, since $rp_A=p_{A}rp_A=p_Ax=x$, $w_j\in A$. As $w_j\in B_H$ we deduce that $\{w_j\}\notin H$. So let $|\mu_j|\geq1$ for all $j$. If $|\alpha_i|=0$ for some $i$ then we let $s_{\g}(\alpha_i)=A_i$. In this case, we have $A\subseteq\cup_is_{\g}(\alpha_i)\cup_js_{\g}(\mu_j)$ (if $v\in A\setminus \cup_is_{\g}(\alpha_i)\cup_js_{\g}(\mu_j)$, then $rp_v=rp_vp_A=p_vx=0$, a contradiction). Thus $A=\big(\cup_is_{\g}(\alpha_i)\cap A\big)\cup\big(\cup_js_{\g}(\mu_j)\cap A\big)$.
Suppose that $\big\{\{v\}: v\in A\big\}\subseteq H$. Since $H$ is hereditary, $A\in H$, which is impossible. Hence there is a vertex $v\in A$ such that $\{v\}\notin H$.

\textbf{Step III}: If $A=\{v\}$, then $\{v\}\in H$.

We go toward a contradiction and  assume $\{v\}\notin H$. Set $v_1:=v$. If $v_1\in B_H$, then there is an edge $e_1\in\gl$ such that $s_{\g}(e_1)=v_1$ and $r_{\g}(e_1)\notin H$. If $v_1\notin B_H$, we have $0<|s^{-1}_{\g}(v)|<\infty$ by Step I. The saturation of $H$ gives an edge $e_1$ with $s_{\g}(e_1)=v_1$ and $r_{\g}(e_1)\notin H$. By Step II, there is a vertex $v_2\in r_{\g}(e_1)$ such that $\{v_2\}\notin H$. We may repeat the argument to choose a path $\gamma=\gamma_1\dots\gamma_L$ for $L=\mathop{\text{max}}_{i,j}\big\{|\beta_i|,|\nu_j|\big\}+1$, such that $s_{\g}(\gamma)=v$ and  $s_{\g}(\gamma_k),r_{\g}(\gamma_k)\notin H$ for  all $1\le k\le L$. Note that $s_{\g}(\gamma_k)\notin A_i$ and so $p_{A_i}p_{s_{\g}(\gamma_k)}=0$ for all $i$, $k$. Moreover, since $r(\gamma_k)\notin H$ we have $p_{w_j}^Hs_{\gamma_k}=0$ for all $j,k$. It follows that $ rs_{\gamma}= (rp_v)s_{\gamma}= xs_{\gamma}=0$,
a contradiction. Therefore $\{v\}\in H$.

\textbf{Step IV}: If $rp_A\in\I$, then $A\in H$.

If $A\notin H$, then by Step II there is a vertex $v\in A$ such that $\{v\}\notin H$, which  contradicts  the Step III. Hence $A\in H$ and consequently $p_A\in\I$, as desired.

Now, we show that $\I$ satisfies assertion (2) of the definition of basic ideal. Note that, by Step IV, we have ${H_{\I}}=H$. Let $w\in B_H$, $r\in R\setminus\{0\}$ and $rp_w^H\in \I$. Using the first part of the lemma, write
$$rp_w^H=x:=\sum\limits_{i=1}^{n}r_is_{\alpha_i}p_{A_i}s_{\beta^*_i} +\sum\limits_{j=1}^{m}s_js_{\mu_j}p^H_{w_j}s_{\nu^*_j},$$
where $A_i\in H$, $w_j\in S$ and $r_i,s_j\in R$ for all $i,j$. Since $rp_w^H=p_w^Hrp_w^H=p_w^Hx=x$ and $p_w^H A_i=0$, We may assume that $|\alpha_i|\geq1$ for all $i$. As $w\in B_H$, we can choose an edge $e\ne\alpha_{1,1},\dots,\alpha_{n,1},\mu_{1,1},\dots,\mu_{m,1}$ such that $s_{\g}(e)=w$ and $r_{\g}(e)\in H$. If $w\notin S$, then $rs_{e^*}=s_{e^*}(rp^H_w)=s_{e^*}x=0 $ which is a contradiction. Therefore $w\in S$ and $p_w^H\in\I$.
\hfill$\Box$
\medskip

\noindent
{\bf Remark 3.2.}
Let $(H,S)$ be an admissible pair in ultragraph $\g$ and let $r\in R\setminus\{0\}$. The argument of Lemma 3.1 implies that
$$H=\big\{A\in \go:rp_A\in I_{(H,S)}\big\}~\text{and}~ S=\big\{w\in B_H:rp_w^H\in I_{(H,S)}\big\}.$$ 

\noindent
{\bf Lemma 3.3.} (cf. \cite[Proposition 3.3]{lar2})
{\it Let $\g$ be an ultragraph and let $R$ be a unital commutative ring. If $(H,S)$ is an admissible pair in $\g$, then $L_R(\g)\cong L_R(\overline{\g})$.}
\vspace{0.05in}\\
\noindent
{\it Proof.}  Suppose that $\{\tilde S,\tilde P\}$ is a universal Leavitt $\overline\g$-family. If we define
\begin{equation}\label{4.1}\begin{array}{ll}
P_A:=\tilde{P}_{\overline A} & \mathrm{for}~~ A\in\go, \\
S_e:=\tilde S_e & \mathrm{for}~~ e\in \gl,\\
S_{e^*}:=\tilde S_{e^*} & \mathrm{for}~~ e\in \gl,
\end{array}
\end{equation}
then it is straightforward to see that $\{S,P\}$ is a Leavitt $\g$-family in $L_R(\overline{\g})$. Note that $L_R(S,P)$ inherits the graded structure of $L_R(\overline\g)$. Since $rP_A\ne0$ for all $A\in\go\setminus\{\emptyset\}$ and $r\in R\setminus\{0\}$, by Corollary 2.15, $L_R(\g)\cong L_R(S,P)$.

We show that $L_R(S,P)=L_R(\overline{\g})$. For  $A\in \go$, $\tilde{P}_{\overline A}=P_A\in L_R(S,P)$. If $v\notin B_H\setminus S$, then $\tilde{P}_{\overline{\{v\}}}=\tilde{P}_v=P_v\in L_R(S,P)$. We note that
\[\overline{s}(e)=
\left\{
\begin{array}{ll}
s_{\mathcal G}(e)' & \mathrm{if}~s_{\mathcal G}(e)\in  B_H\setminus S ~\mathrm{and}~r_{\mathcal G}(e)\in H,\\
s_{\mathcal G}(e) & \text{otherwise},\\
\end{array}
\right.
\]
for every $e\in{\g}^1$. Thus for $w\in B_H\setminus S$, we have $0<|\overline s^{-1}(w)|<\infty$ and $r_{\g}(e)\notin H$ for $e\in \overline s^{-1}(w)$. Hence
$$\tilde{P}_w=\sum\limits_{\overline s(e)=w}\tilde S_e\tilde S_{e^*}=\sum\limits_{ s_{\g}(e)=w,\,\,  r_{\g}(e)\notin H} S_e S_{e^*}\in L_R(S,P).$$
Also,
$$\tilde P_{w'}=\tilde P_{\{w,w'\}}-\tilde P_w=P_w-\sum\limits_{ s_{\g}(e)=w,\,\,  r_{\g}(e)\notin H} S_e S_{e^*}=P_w^H\in L_R(S,P).$$
Thus,  by Remark 2.6, $\tilde P_{A}\in  L_R(S,P)$ for all $A\in\overline\g$. Since $\tilde S_e,\tilde S_{e^*}\in L_R(S,P)$, we deduce that $L_R(S,P)=L_R(\overline{\g})$. Consequently, $L_R(\g)\cong L_R(\overline{\g})$.
\hfill$\Box$
\medskip

\noindent
{\bf Theorem 3.4.} (cf. \cite[Theorem 3.10]{lar1})
{\it Let $\g$ be an ultragraph  and let $R$ be a unital commutative ring. Then
	\begin{enumerate}
		\item For any  admissible pair $(H,S)$ in $\g$, we have $\lqg\cong L_R(\g)/\I$.
		\item  The map $(H,S)\mapsto\I$ is a bijection from the set of all admissible pairs of $\g$  to the set of all graded basic ideals of $L_R(\g)$.
\end{enumerate}}
\noindent
{\it Proof.} (1) Let $\{\tilde S,\tilde P\}$ be a universal Leavitt $\overline\g$-family and let $L_R(\g)=L_R(S,P)$, where $\{S,P\}$ is the Leavitt ${\g}$-family as defined in  \ref{4.1}.  Define
\begin{equation}\label{4.2}\begin{array}{ll}
Q_{[A]}:=\tilde P_A+I_{(H,S)} & \mathrm{for}~~A\in\overline{\g}^0, \\
T_{e}:=\tilde S_e+I_{(H,S)} &  \mathrm{for}~~e\in\qgl, \\
T_{e^*}:=\tilde S_{e^*}+I_{(H,S)} & \mathrm{for}~~ e\in\qgl.
\end{array}
\end{equation} 
It can be shown that the family $\{Q_{[A]},T_e,T_{e^*}:[A]\in \qgo,e\in \qgl\}$ is a Leavitt  $\qg$-family that generates $L_R(\overline{\g})/\I$. Furthermore, by Remark 3.2, $rQ_{[A]}\ne0$ for all $[A]\in\qgo\setminus\{[\emptyset]\}$ and $r\in R\setminus \{0\}$.

Now we use the universal property of $\lqg$ to get an $R$-homomorphism $\pi:\lqg\rightarrow L_R(\overline\g)/\I$ such that $\pi(t_e)=T_e$, $\pi(t_{e^*})=T_{e^*}$ and $\pi(rq_{[A]})=rQ_{[A]}\ne0$ for $[A]\in\qgo\setminus\{[\emptyset]\}$, $e\in\qgl$ and $r\in R\setminus \{0\}$. Since $\I$  is a graded ideal, the quotient $L_R(\overline\g)/\I$ is graded. Moreover, the elements  $Q_{[A]}$, $T_e$ and $T_{e^*}$ have degrees 0, 1 and -1 in $L_R(\overline\g)/\I$, respectively and thus $\pi$ is a graded homomorphism. It follows from the graded-uniqueness Theorem 2.14 that $\pi$ is injective. Since $L_R(\overline\g)/\I$ is generated by $\{\ Q_{[A]},T_{e},T_{e^*}:[A]\in\qgo,e\in\qgl\}$, we deduce that $\pi$ is also surjective. Hence $\lqg\cong L_R(\overline\g)/\I\cong L_R(\g)/\I$.

(2) The injectivity of the map $(H,S)\mapsto \I$ is a consequence of Remark 3.2. To see that it is onto, let $I$ be a graded basic ideal in $L_R(\g)$. Then $I_{(H_I,S_I)}\subseteq I$. Consider the ultragraph $\overline\g$ with respect to admissible pair $(H_I,S_I)$.
Since $I$ is a graded ideal, the quotient ring $L_R(\overline\g)/I$ is graded. Let $\pi: L_R(\g/(H_I,S_I))\cong L_R(\overline\g)/I_{(H_I,S_I)}\rightarrow L_R(\overline\g)/I$ be the quotient map. For $(H_I,S_I)$, consider $\{\tilde S,\tilde P\}$ and $\{T,Q\}$ as defined in Equations \ref{4.1} and \ref{4.2}, respectively. Since $I$ is basic, we have $rp_A\notin I$ and $rp_w^{H_I}\notin I$ for $A\in \go\setminus H_I$, $w\in B_{H_I}\setminus S_{I}$ and $r\in R\setminus\{0\}$.

We show that $\pi(rq_{[A]})=\pi(rQ_{[A]})=r\tilde P_{A}+I\ne0$ for all $[A]\in\qgo\setminus\{[\emptyset]\}$ and $r\in R\setminus\{0\}$. Fix $r\in R\setminus\{0\}$. We know that $\qgo=\big\{[A]:A\in\overline{\g}^0\big\}$. Let $A=\overline B$ for some $B\in\go$. If $[A]\ne[\emptyset]$, then $B\notin H_I$ and thus  $r\tilde P_{A}=rp_B\notin I$. Therefore $\pi(rq_{[A]})\ne0$. If $w\in B_{H_I}\setminus S_I$, then $r\tilde P_{w'}=rp_w^{H_I}\notin I$. Hence $\pi(rq_{[w']})\ne0$. In view of Remark 2.10, we deduce that $\pi(rq_{[A]})\ne0$ for every $[A]\in\qgo\setminus\{[\emptyset]\}$. Furthermore, $\pi$ is a graded homomorphism. It follows from  Theorem 2.14 that the quotient map $\pi$ is injective. Hence $I=I_{(H_I,S_I)}$.
\hfill$\Box$
\medskip

As we have seen in the proof of Theorem 3.4, if $(H,S)$ is an admissible pair in  $\g$, then there is  a Leavitt  $\qg$-family $\{T,Q\}$ in $L_R(\g)/I_{(H,S)}$ such that $rQ_{[A]}\ne0$ for all $[A]\in\qgo\setminus\{[\emptyset]\}$ and $r\in R\setminus\{0\}$. So we have the following proposition.
\vspace{0.1in}\\
\noindent
{\bf Proposition 3.5.} Let $(H,S)$ be an admissible pair in ultragraph $\g$ and let $R$ be a unital commutative ring. If $\{t,q\}$ is the universal Leavitt  $\qg$-family, then  $rq_{[A]}\ne0$ for every $[A]\in\qgo\setminus\{[\emptyset]\}$ and every $r\in R\setminus\{0\}$.

%%%%%%%%%%%%%%%%%%%%%%%%%%%%%%%%%%%%%%%
\section{Condition (K)}\label{s4}

In this section we recall Condition (K) for ultragraphs and 
we consider the ultragraph $\g$ that satisfy Condition (K) to describe that  
all basic ideals of $L_R(\g)$ are graded.

Let $\g=(G^0,\gl,r_{\g},s_{\g})$ be an ultragraph and let $v\in G^0$. A \emph{first-return path based at} $v$ in $\g$ is a loop $\alpha=e_1e_2\cdots e_n$ such that $s_{\g}(\alpha)=v$ and $s_{\g}(e_i)\ne v$ for $i=2,3,\dots,n$.
\vspace{0.1in}\\
\noindent
{\bf Definition 4.1.} \cite[Section 7]{kat2} 	An ultragraph $\g$ satisfies Condition (K) if every vertex in $G^0$ is either the base of no first-return path or it is the base of at least two first-return paths.
\medskip

Let $\qg$ be a quotient ultragraph. By rewriting Definition 2.2 for $\qg$, one can define the hereditary property for the subcollections of $\qgo$. More precisely, a subcollection $K\subseteq\qgo$ is called \emph{hereditary} if satisfying the following conditions:
\begin{enumerate}
	\item $s(e)\in K$ implies $r(e)\in K$ for all $e\in\qgl$.
	\item $[A]\cup [B]\in K$ for all $[A],[B]\in K$.
	\item $[A]\in K$, $[B]\in\qgo$ and $[B]\subseteq [A]$, imply $[B]\in K$.
\end{enumerate}
For any hereditary subcollection $K\subseteq\qgo$, the ideal  $I_K$  in $\qg$ generated by $\{q_{[A]}:[A]\in K\}$ is equal to $$\text{span}_R\big\{t_{\alpha}q_{[A]}t_{\beta^*}\in\lqg:\alpha,\beta\in\text{Path}(\qg)~\text{and}~[A]\in {K}\big\}.$$

\noindent
{\bf Lemma 4.2.}
{\it Let $\qg$ be a quotient ultragraph and let $R$ be a unital commutative ring. If $\qg$ does not satisfy Condition (L), then $\lqg$ contains a non-graded ideal $I$ such that $rq_{[A]}\notin I$ for all $[A]\in\qgo\setminus\{[\emptyset]\}$ and  $r\in R\setminus\{0\}$.}
\vspace{0.05in}\\
\noindent
{\it Proof.} Suppose that $\qg$ contains a closed path $\gamma:=e_1e_2\cdots e_n$ without exits and $r(e_i)=s(e_{i+1})$ for $1\le i\le n$ where $e_{n+1}:=e_1$. Thus we can assume that $s(e_i)\ne s(e_j)$ for all $i,j$. Let $s(e_i)=[v_i]$ for $1\le i\le n$ and let $X$ be the subalgebra of $L_R\big(\qg\big)$ generated by $\{t_{e_i},t_{e_i^*}, q_{[v_i]}:1\le i\le n\}$.

{\bf Claim 1:} The subalgebra $X$ is isomorphic to the Leavitt path algebra $L_R(E)$, where $E$ is the graph containing a single simple closed path of length $n$, that is,
$$E^0=\{w_1,\dots,w_n\},~~~~~ E^1=\{f_1,\dots,f_n\},$$
$$r(f_i)=s(f_{i+1})=w_{i+1}~~ \text{for}~~1\le i\le n~~\text{where}~~w_{n+1}:=w_1.$$

Proof of Claim 1. Set $P_{w_i}:=q_{[v_i]}$, $S_{f_i}:=t_{e_i}$ and $S_{f_i^*}:=t_{e_i^*}$ for all $i$. Then $\{S,P\}$ is a Leavitt  $E$-family in $X$ such that $L_R(S,P)=X$. So there is a homomorphism $\pi:L_R(E)\rightarrow X$ such that $\pi(s_{f_i})=t_{e_i}$, $\pi(s_{f_i^*})=t_{e_i^*}$ and $\pi(p_{w_i})=q_{[v_i]}$ for all $i$. Since $\pi$ preserves the degree of generators, $\pi$ is a graded homomorphism. By Proposition 3.5, $\pi(rp_{w_i})=rq_{[v_i]}\ne0$ for all $r\in R\setminus\{0\}$ and all $i$. Now apply \cite[Theorem 5.3]{tom2} to obtain the injectivity of $\pi$. Hence, $X\cong L_R(E)$ which proves the  claim.

By \cite[Lemma 7.16]{tom2},  $L_R(E)$ contains an ideal $J$ that is basic but not graded. Thus $\pi(J)$ is a non-graded ideal of $X$.  

{\bf Claim 2:} $rq_{[A]}\notin\pi(J)$ for every  $r\in R\setminus\{0\}$ and  $[A]\in \qgo\setminus[\emptyset]$.

Proof of Claim 2. We show that $rp_{w_i}\notin J$ for all $r\in R\setminus\{0\}$ and all  $i$. Let $rp_{w_j}\in J$ for some $r\in R\setminus\{0\}$ and some $1\leq j\leq n$. Since $J$ is a basic ideal, $p_{w_j}\in J$. For every $1\leq i\leq n$, there is a path $\alpha\in\mathrm{Path}(E)$ such that $s(\alpha)=w_j$ and $r(\alpha)=w_i$. Hence $p_{w_i}=s_{\alpha_i^*}p_{w_j}s_{\alpha_i}\in J$ for all $i$. Thus $J=L_R(E)$, contradicting the hypothesis that $J$ is non-graded.

Suppose that there exist $r\in R\setminus\{0\}$ and  $[A]\in \qgo\setminus[\emptyset]$ such that $rq_{[A]}\in\pi(J)$. If we identify $X$ with the subalgebra
$${\text{span}_R}\big\{t_{\alpha}t_{\beta^*}\in\lqg:s(\alpha),s(\beta)\in\{[v_1],\dots,[v_n]\}\big\}$$
of $L_R\big(\qg\big)$, then $rq_{[A]}=r_{1}t_{\alpha_1}t_{\beta_1^*}+\ldots+r_mt_{\alpha_m}t_{\beta_m^*}$, where $r_k\in R\setminus\{0\}$ and $s(\alpha_k),s(\beta_k)\in\{[v_1],\dots,[v_n]\}$. Thus 
$$rq_{[A]}=rq_{[A]}q_{[A]}=r_1q_{[A]}t_{\alpha_1}t_{\beta_1^*}+\ldots+r_mq_{[A]}t_{\alpha_m}t_{\beta_m^*}\ne0,$$
which implies that $[A]\cap(\cup_{i=1}^n[v_i])\ne[\emptyset]$. So there is $j$ such that $[v_j]\subseteq[A]$. In this case we  have $rq_{[v_j]}=q_{[v_j]}rq_{[A]}\in\pi(J)$. Hence $rp_{w_j}=\pi^{-1}(rq_{[v_j]})\in J$, a contradiction.

Let
$$K:=\Big\{\bigcup\limits_{k=1}^m[v_{n_k}]:1\le m\le n~\text{and}~n_k\in\{1,2,\dots,n\}\Big\}$$
and set $[A]:=\cup_{i=1}^n[v_i]$. Since $\gamma$ is a closed path without exits and $r(e_i)=s(e_{i+1})$, $K$ is a hereditary subset of $\qgo$. 

{\bf Claim 3:} If $rq_{[v]}\in I_K$ for some $[v]\in \qgo$ and some $r\in R\setminus\{0\}$, then there exists $\alpha\in\text{Path}\big(\qg\big)$ such that $s(\alpha)=[v]$ and $r(\alpha)\cap A\ne[\emptyset]$.

Proof of Claim 3. Suppose that
$$rq_{[v]}=x:=r_1t_{\alpha_1}q_{[v_{n_1}]}t_{\beta^*_1}+\ldots +r_mt_{\alpha_m}q_{[v_{n_m}]}t_{\beta^*_m},$$
where  $r_i\in R\setminus\{0\}$ and $\alpha_i,\beta_i\in\text{Path}\big(\qg\big)$. If there is no such path, then we can choose an edge $e$ such that $s(e)=[v]$ and $r(e)\cap[A]=[\emptyset]$. Since $rq_{r(e)}\in I_K$, there is a vertex $[w]\subseteq r(e)$ such that $rq_{[w]}\in I_K$. By a continuing process, one can  choose a path $\eta=\eta_1\cdots\eta_L$ for $L=\mathop{\text{max}}_{i}\big\{|\beta_i|\big\}+1$, such that $s(\eta)=[v]$ and $r(\eta_k)\cap [A]=[\emptyset]$ for all $1\le k\le L$. Thus we have $rt_{\eta}= (rq_{[v]})t_{\eta}= xt_{\eta}=0,$ a contradiction. Therefore, there is a path $\alpha$ such that $s(\alpha)=[v]$ and $r(\alpha)\cap [A]\ne[\emptyset]$.

{\bf Claim 4:} $I:=I_K\pi(J)I_K$ is a non-graded ideal of $\lqg$.

Proof of Claim 4. Suppose that $I$ is graded. Thus $I=\bigoplus I_n$, where $I_n=I\cap\lqg_n$. For $x\in \pi(J)$ and $y\in\lqg$ we have $q_{[A]}xq_{[A]}=x$ and $q_{[A]}y$,$yq_{[A]}\in X$. Thus $\pi(J)\subseteq I$. Let $x\in \pi(J)$. Since $I=\bigoplus I_n$ we have
$x=r_{n_1}y_{n_1}x_{n_1}z_{n_1}+\ldots+ r_{n_m}y_{n_m}x_{n_m}z_{n_m}$, where $x_{n_k}\in\pi(J)$, $y_{n_k},z_{n_k}\in I_K$ and $y_{n_k}x_{n_k}z_{n_k}\in I_{n_k}$. Thus
$$x=q_{[A]}xq_{[A]}=r_{n_1}q_{[A]}y_{n_1}x_{n_1}z_{n_1}q_{[A]}+\ldots+ r_{n_m}q_{[A]}y_{n_m}x_{n_m}z_{n_m}q_{[A]}.$$
Since $q_{[A]}y_{n_k}$,$z_{n_k}q_{[A]}\in X$, we have $r_{n_k}q_{[A]}y_{n_k}x_{n_k}z_{n_k}q_{[A]}\in\pi(J)\cap X_n$, where $X_n:=X\cap\lqg_n$. Hence $\pi(J)$ is a graded ideal of $X$, a contradiction.

Finally, we show that $rq_{[B]}\notin I$ for all $[B]\in\qgo\setminus\{[\emptyset]\}$ and  $r\in R\setminus\{0\}$. Assume $rq_{[B]}\in I$ for some $r\in R\setminus\{0\}$ and some $[B]\in\qgo\setminus\{[\emptyset]\}$. Since $rq_{[B]}\in I_K$, there is a vertex $[v]\subseteq [B]$ such that $rq_{[v]}\in I\subseteq I_K$. By Claim 2, there is a path $\alpha$ such that $s(\alpha)=[v]$ and $[v_i]\subseteq r(\alpha)\cap[A]$ for some $i\in\{1,\ldots,n\}$. So $rq_{[v_i]}=t_{\alpha^*}rq_{[v]}t_{\alpha}q_{[v_i]}\in I$. If we consider $rq_{[v_i]}$ in terms of its representation in $I_K\pi(J)I_K$, then one can show that $rq_{[v_i]}\in X\pi(J)X=\pi(J)$ which is a contradiction. Therefore, we have $rq_{[B]}\notin I$ for every $r\in R\setminus\{0\}$ and every $[B]\in\qgo\setminus\{[\emptyset]\}$.
\hfill$\Box$
\medskip

\noindent
{\bf Theorem 4.3.} (cf. \cite[Theorem 3.18]{lar1}, \cite[Theorem 7.17]{tom2})
{\it 	Let $\g$ be an ultragraph and let R be a unital commutative ring. Then $\g$ satisfies Condition (K) if and only if every basic ideal in $L_R(\g)$ is graded.}
\vspace{0.05in}\\
\noindent
{\it Proof.} 	Suppose that $\g$ satisfies Condition (K). If $I$ is a basic ideal of $L_R(\g)$, then by Theorem 3.4, we have the quotient map
$$\pi: L_R(\g/(H_I,S_I))\cong L_R(\g)/I_{(H_I,S_I)}\rightarrow L_R(\g)/I,$$
such that $\pi(rq_{[A]})\ne0$ for all $[\emptyset]\ne[A]\in \qgo$ and all $r\in R\setminus\{0\}$. By Proposition \cite[Proposition 6.2]{lar2} and the Cuntz-Krieger uniqueness Theorem 2.17, $\pi$ is injective and $I=I_{(H_I,S_I)}$. It follows from Lemma 3.1 that $I$ is a graded ideal.

Conversely, if $\g$ does not satisfy Condition (K), then \cite[Proposition 6.2]{lar2} and Lemma 4.2 imply that there exist $H$ and $S$ such that $L_R(\qg)\cong L_R(\g)/I_{(H,S)}$ contains a non-graded ideal $J$ such that $rq_{[A]}\notin J$ for all $[A]\in\qgo\setminus\{[\emptyset]\}$ and all $r\in R\setminus\{0\}$. If $\pi:L_R(\g)\rightarrow  L_R(\qg)$ is the quotient map, then $I:=\pi^{-1}(J)$ is a non-graded ideal of $L_R(\g)$. Similar to the last paragraph  of the  proof of Theorem 3.4, we have  $rp_{A}\notin I$   for  $\emptyset\ne A\notin H$ and  $r\in R\setminus\{0\}$, that is, $H=H_{I}:=\{A\in\go:p_A\in I\}$. Also, $rp^H_w\notin I$ for $w\in B_H\setminus S$ and $r\in R\setminus\{0\}$. Consequently $I$ is a basic ideal. Therefore $\g$ satisfies Condition (K).
\hfill$\Box$
\medskip

\noindent
{\bf Corollary 4.4.} Let $\g$ be an ultragraph and let $R$ be a unital commutative ring. If $\g$ satisfies Condition (K), then the map $(H,S)\mapsto\I$ is a bijection from the set of all admissible pairs of $\g$ into the set of all basic ideals of $L_R(\g)$.

%%%%%%%%%%%%%%%%%%%%%%%%%%%%%%%%%%%%%%%%

\section{Exel-Laca $R$-algebras}\label{s5}

The  Exel-Laca algebras $\mathcal O_A$ are generated by partial isometries whose relations are determined by a countable $\{0,1\}$-valued  matrix $A$ with no identically zero rows  \cite[Definition 8.1]{exe}. In this section, we define an algebraic version of the Exel-Laca algebras.
\vspace{0.1in}\\
\noindent
{\bf Definition 5.1.} Let $I$ be a countable set and let $A$ be an $I\times I$ matrix with entries in $\{0, 1\}$. The ultragraph $\g_A:=(G^0_A,\gl_A,r,s)$ is defined by $G^0_A:=\{v_i:i\in I\}$, $\gl_A:=I$, $s(i):=v_i$ and $r(i):=\{v_j:A(i,j)=1\}$ for all $i\in I$.
\medskip

By \cite[Theorem 4.5]{tom1}, the  Exel-Laca algebra $\mathcal O_A$ is canonically isomorphic to $C^*(\g_A)$.  Motivated by this fact we give the following definition.
\vspace{0.1in}\\
\noindent
{\bf Definition 5.2.} Let $R$ be a unital commutative ring, let $I$ be a countable set and let $A$ be a $\{0,1\}$-valued matrix over $I$ with no identically zero rows. The  \emph{Exel-Laca $R$-algebra  associated} to $A$, denoted by $\el$, is defined by $\el:=L_R(\g_A)$.
\medskip

\noindent
{\bf Example 5.3} Let $R$ be a unital commutative ring and let $A=\big(A(i,j)\big)$ be an $n\times n$ matrix, with with no identically zero rows and $A(i,j)\in\{0_R,1_R\}$ for all $i,j$. The \emph{Cuntz-Krieger $R$-algebra associated} to $A$, as defined in \cite[Example 2.5]{ara}, is the free $R$-algebra generated by a set $\{x_1,y_1,\dots,x_n,y_n\}$, modulo the ideal generated by the following relations:
\begin{enumerate}
	\item $x_iy_ix_i=x_i$ and $y_ix_iy_i=y_i$ for all $i$,
	\item $y_ix_j=0$ for all $i\ne j$,
	\item $y_ix_i=\sum_{j=1}^nA(i,j)x_jy_j$ for all $i$,
	\item $\sum_{j=1}^nx_jy_j=1$.
\end{enumerate}
The Cuntz-Krieger $R$-algebra associated to $A$ is denoted by $\mathcal{CK}_A(R)$. We will show that $\mathcal{CK}_A(R)\cong\el$.  Note that $\g_A$ is a finite ultragraph. If we define $P_B:=\sum_{v_i\in B}x_iy_i$, $S_i:=x_i$ and $S_{i^*}:=y_i$ for all $i$ and all $B\in\g_A^0$, then it is easy to verify that $\{P_B,S_i,S_{i^*}: B\in\g_A^0, 1\leq i\leq n\}$ is a Leavitt $\g_A$-family.

On the other hand, if $\{s,p\}$ is the universal Leavitt $\g_A$-family for $L_R(\g_A)$, then the elements $X_i:=s_i$ and $Y_i:=s_{i^*}$, $1\leq i\leq n$, satisfy relations (1)-(4). Now the universality of $\mathcal{CK}_A(R)$ and $L_R(\g_A)$ conclude that  $\mathcal{CK}_A(R)\cong\el$.
\medskip

Note that for any $\{0,1\}$-valued matrix $A$ with no identically zero rows, the ultragraph $\g_A$ contains no singular vertices. So, each Exel-Laca algebra is the Leavitt path algebra of a non-singular ultragraph. In Theorem 5.7 below, we will prove the converse. For this, given any ultragraph $\g=(G^0,\gl,r,s)$ with no singular vertices, we first reform $\g$ to an ultragraph $\widetilde{\g}$ as follows. Associated to each $v\in G^0$ and $e\in s^{-1}(v)$, we introduce a vertex $v_e$ and set $\widetilde A:=\{v_e:v\in A~\text{and}~e\in s^{-1}(v)\}$ for every $A\in\go$. Next we  define the ultragraph $\widetilde\g:=(\widetilde G^0,\widetilde{\mathcal G}^1,\widetilde r,\widetilde s)$, where $\widetilde G^0:=\big\{v_e:v\in G^0 ~\text{and}~e\in s^{-1}(v)\big\}$, $\widetilde{\mathcal G}^1:=\gl$, the source map $\widetilde s(e):=s(e)_e$ and {\tiny } the range map $\widetilde r(e):=\widetilde{r(e)}$, where $\widetilde{r(e)}=\{v_e:v\in r(e)~\text{and}~e\in s^{-1}(v)\}$.
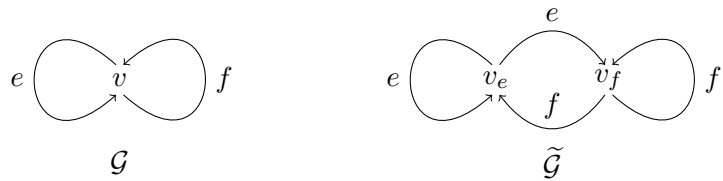
\begin{center}
	\begin{tikzpicture}
	\draw (-1,0) node(x)  {$v$} (-1,0);
	\draw[->] (-1.05,.2) .. controls (-2.5,1.5) and (-2.5,-1.5) ..node[left] {$e$}(-1.05,-.2);
	\draw[->] (-.95,-.2) .. controls (.5,-1.5) and (.5,1.5) ..node[right] {$f$}(-.95,.2);
	\draw (-1,-1.1) node(x)  {$\g$} (-1,-1.1);
	
	\draw (4,0) node(x)  {$v_e$} (4,0);
	\draw[->] (3.95,.2) .. controls (2.5,1.5) and (2.5,-1.5) ..node[left] {$e$}(3.95,-.2);
	\draw[->] (4.05,.2) .. controls (4.5,.8) and (5,.8) ..node[above] {$e$}(5.45,.2);
	\draw (5.5,0) node(x)  {$v_f$} (5.5,0);
	\draw[->] (5.55,-.2) .. controls (7,-1.5) and (7,1.5) ..node[right] {$f$}(5.55,.2);
	\draw[->] (5.45,-.2) .. controls (5,-.8) and (4.5,-.8) ..node[above] {$f$}(4.05,-.2);
	\draw (4.75,-1.1) node(x)  {$\widetilde\g$} (4.75,-1.1);
	\end{tikzpicture}
\end{center}

\noindent
{\bf Remark 5.4.} 	We notice that each vertex $v_e$ in $\widetilde{\g}$ emits exactly one edge $e$. Moreover, Lemma 2.5 implies that
$$\widetilde{\mathcal G}^0=\{\widetilde A\cup F:A\in\go~\text{and}~F~\text{is}~\text{a}~\text{finite}~\text{subset}~\text{of}~\widetilde G^0\}.$$

\noindent
{\bf Lemma 5.5.}
{\it 	Let $\g$ be a row-finite ultragraph with no singular vertices and let $R$ be a unital commutative ring. Then $L_R(\widetilde\g)$ is isomorphic to  $L_R(\g)$ as algebras.}
\vspace{0.05in}\\
\noindent
{\it Proof.} Let $\{t,q\}$ be the universal Leavitt  $\widetilde\g$-family for $L_R(\widetilde\g)$. Since $\g$ is row-finite, for every $v\in G^0$ we have $|\{v_e:e\in\gl\}|<\infty$. Define

\begin{equation*}\begin{array}{ll}
P_{A}:=q_{\widetilde A} & \mathrm{for}~~ A\in \go,\\
S_e:=t_e & \mathrm{for}~~ e\in \gl,\\
S_{e^*}:=t_{e^*} & \mathrm{for}~~ e\in \gl.
\end{array}
\end{equation*}
By Remark 5.4, for each $A\in\go$,  the idempotent $P_A$ satisfies the condition (1) of Definition 2.1.
It is straightforward to check that $\{S,P\}$ is a Leavitt $\g$-family in $L_R(\widetilde\g)$. For example, to verify condition (2) of Definition 2.1 suppose that $e\in\gl$. Then $S_eP_{r(e)}=t_{e}q_{\widetilde r(e)}=t_e=S_e$ and
$$P_{s(e)}S_e=\big(\sum_{f\in\widetilde\g^1}q_{s(e)_f}\big)q_{\widetilde s(e)}t_e=\big(\sum_{f\in\widetilde\g^1}q_{s(e)_f}\big)q_{s(e)_e}t_e=q_{s(e)_e}t_e =t_e=S_e.$$

Thus there is an $R$-algebra homomorphism $\phi:L_R(\g)\rightarrow L_R(\widetilde\g)$ such that $\phi(s_e)=t_e$, $\phi(s_{e^*})=t_{e^*}$ and $\phi(rp_v)=r\sum_{\substack{e\in\widetilde\g^1}}q_{v_e}\ne0$ for $v\in G^0$, $e\in \gl$ and $r\in R\setminus \{0\}$. So $r\phi(p_A)\ne0$ for all $A\in\go\setminus\{\emptyset\}$. As $\phi$ is a graded homomorphism, Corollary 2.15 implies that $\phi$ is injective. Moreover, for any $v_e\in\widetilde G^0$, we have $(\widetilde s)^{-1}(v_e)=\{e\}$ and hence $S_eS_{e^*}=t_et_{e^*}=q_{v_e}$. Therefore, $\phi$ is an isomorphism and we get the result.
\hfill$\Box$
\medskip

\noindent
{\bf Definition 5.6.} \cite[Definition 2.4]{tom1} 	Let $\g$ be an ultragraph. The edge matrix of an ultragraph $\g$ is the $\gl\times\gl$ matrix $A_\g$ given by
\begin{equation*}
A_\g(e,f):=
\begin{cases}
1 \,\,\,\,\,\,\,\,\mathrm{if}\,\, s(f)\in r(e),\\
0\,\,\,\,\,\,\,\, \text{otherwise}.
\end{cases}
\end{equation*}
\medskip

We can check that, if $\g$ is an ultragraph with no singular vertices, then $\g_{A_\g}=\widetilde{\g}$. So, by Lemma 5.5 we have the following.
\vspace{0.1in}\\
\noindent
{\bf Theorem 5.7.}
{\it Let $R$ be a unital commutative ring and let $\g$ be an ultragraph with no singular vertices. Then $L_R(\g)$ is isomorphic to $\mathcal{EL}_{A_\g}(R)$.}
\medskip

\noindent
{\bf Example 5.8.}
Suppose that $E$ is a graph with one vertex $E^0=\{v\}$ and $n$ edges $E^1=\{e_1,\dots.e_n\}$. We know that $L_R(E)$ is isomorphic to the \emph{Leavitt algebra} $L_R(1,n)$. If  all  entries of matrix $A_{n\times n}$ is 1, then $A=A_E$. Hence, Theorem 5.7 shows that $\el\cong L_R(1,n)$.
\medskip

\noindent
{\bf Remark 5.9.}
	Theorem 5.7 shows that the Leavitt path algebras of ultragraphs with no singular vertices are precisely the Exel-Laca $R$-algebras.
\medskip

\noindent
{\bf Remark 5.10.} 
We know that for every directed graph $E$, $L_{\mathbb{C}}(E)$ is isomorphic to a dense $*$-subalgebra of the $C^*$-algebra $C^*(E)$ (see  \cite[Section 7]{tom4}).
Using a similar argument to \cite[Theorem 7.3]{tom4}, one can  show that $L_{\mathbb{C}}(\g)$ is isomorphic to a dense $*$-subalgebra of the $C^*$-algebra $C^*(\g)$ introduced in \cite{tom1}. In particular, for any countable $\{0,1\}$-valued matrix $A$, the Exel-Laca algebra $\mathcal{EL}_A(\mathbb C)$ is isomorphic to a dense $*$-subalgebra of $\mathcal O_A$ \cite{exe}.
\medskip

Last part of this section is an  example to emphasize that the class of the Leavitt path algebras  of ultragraphs is strictly larger than the class of the Leavitt path algebras of directed graphs as well as  the class of algebraic Exel-Laca algebras.
\vspace{0.1in}\\
\noindent
{\bf Example 5.11.}
Let $R=\mathbb{Z}_2$ and let $\g$ be the ultragraph 
$$
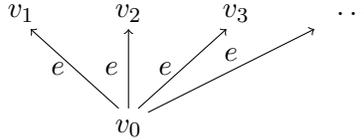
\begin{tikzpicture} [x=37pt,y=30pt]
\draw [->] (-.1,.25) -- node[ left] {$e$} (-1,1.25);
\draw (-1.1,1.45) node(x)  {$v_1$} (-1.1,1.45);
\draw [->] (0,.25) -- node[ left] {$e$} (0,1.25);
\draw (0,1.45) node(x)  {$v_2$} (0,1.45);
\draw [->] (.1,.25) -- node[ left] {$e$} (1,1.25);
\draw (1.1,1.45) node(x)  {$v_3$} (1.1,1.45);
\draw [->] (.2,.20) -- node[above ] {$e$} (1.9,1.25);
\draw (2.3,1.45) node(x)  {$\ldots$} (2.3,1.45);
\draw (0,0) node(x)  {$v_0$} (0,0);
\end{tikzpicture}
$$
with one ultraedge $e$ such that $s(e)=v_0$ and $r(e)=\{v_1,v_2,\ldots\}$. Note that $\g$ satisfies Condition (K) and so all ideals of $L_{\mathbb Z_2}(\g)$ are basic and graded by Theorem 4.3. Moreover, $L_{\mathbb Z_2}(\g)$ is unital, because by Lemma 2.12 $G^0=r(e)\cup \{v_0\}\in \go$. We claim that the algebra $L_{\mathbb Z_2}(\g)$ is not isomorphic to the Leavitt path algebra  of a graph. Suppose on the contrary that there is a such graph $E$.  Let $\{t,q\}$ be the canonical generators for $L_{\mathbb Z_2}(E)$. If $E$  has a loop $\alpha$, then for $t_\alpha\in L_{\mathbb Z_2}(E)$ we have $t_{\alpha}^n\ne t_{\alpha}^m$ for all $m\ne n$, but one can show that $L_{\mathbb Z_2}(\g)$ does not include such member. Thus $L_{\mathbb Z_2}(E)$ does not have any loop and consequently $E$ satisfies Condition (K) for directed graphs. Hence by \cite[Theorem 3.18]{lar1} any ideal of $L_{\mathbb Z_2}(E)$ is graded. Since $L_{\mathbb Z_2}(E)$ is unital, $E^0$ must be finite and hence $L_{\mathbb Z_2}(E)$ contains finitely many (graded) ideals  by Theorem 3.4(2),  which contradicts the fact that $L_{\mathbb Z_2}(\g)$ has infinitely many pairwise orthogonal ideals $I_{\{v_n\}}$ for $n\geq 1$.

Now assume $L_{\mathbb Z_2}(\g)\cong\mathcal{EL}_{A}(\mathbb Z_2)$ for some matrix $A$ and let $\phi:L_{\mathbb Z_2}(\g)\rightarrow \mathcal{EL}_{A}(\mathbb Z_2)$ be an isomorphism. We may consider the ideal $I_{\{v_1\}}$ in $L_{\mathbb Z_2}(\g)$. Since $v_1$ is a sink, we have
$$I_{\{v_1\}}=\text{span}_{\mathbb Z_2}\{p_{v_1},s_ep_{v_1},p_{v_1}s_{e^*},s_{e}p_{v_1}s_{e^*}\},$$
and so $|I_{\{v_1\}}|<\infty$. On the other hand, we know that $\g_A$ has no sinks. Hence the graded ideal $\phi(I_{\{v_1\}})$ in $L_{\mathbb Z_2}(\g_A)$ includes infinitely many elements, which is a contradiction.
\medskip
\vspace{0.2in}\\
{\bf Acknowledgment.} The authors are  grateful to the referee for carefully reading the paper, pointing out a number of misprints and for some helpful comments.

\footnotesize{

}
\end{document}